\documentclass[11pt]{article}
\usepackage{amssymb}
\usepackage{amsmath}
\usepackage{amsthm}
\usepackage{geometry} 
\usepackage{hyperref}
\usepackage{parskip}
\usepackage{comment}
\usepackage{enumerate}
\usepackage[toc]{appendix}
\usepackage[scr=rsfso,frak=euler]{mathalfa}

\usepackage{soul,xcolor}
\usepackage{soul,xcolor}

\usepackage{chemformula}
 
\usepackage{lineno}

\usepackage{enumerate}

\usepackage{bbm}

\RequirePackage{geometry}
\theoremstyle{remark}
\newtheorem{remark}{Remark}

\usepackage[
    backend=biber,
    style=numeric,
    sorting=nyt,
    maxbibnames=10,
  ]{biblatex}
  \renewbibmacro{in:}{%
  \ifentrytype{article}
    {}
    {\bibstring{in}%
     \printunit{\intitlepunct}}}
\usepackage{authblk}

\begin{document}  

\title{Path-dependent McKean PDEs with reaction: a discussion on probabilistic interpretations and particle approximations
}

\author{Daniela Morale\thanks{Dept. of Mathematics, University of Milano, \href{mailto:daniela.morale@unimi.it}{daniela.morale@unimi.it}} , Leonardo Tarquini\thanks{Dept. of Mathematics, University of Oslo, \href{mailto:leonardt@math.uio.no}{leonardt@math.uio.no}} , Stefania Ugolini\thanks{Dept. of Mathematics, University of Milano, \href{mailto:stefania.ugolini@unimi.it}{stefania.ugolini@unimi.it}}}

\maketitle

\begin{abstract}
    In this paper, we discuss and compare two probabilistic approaches for associating a stochastic differential equation with a McKean-type partial differential equation featuring a reaction term and path-dependent coefficients given in \cite{arxivTarquiniUgolini1,2025_morale_tarquini_ugolini}. The non-conservative nature of the macroscopic dynamics leads to two possible interpretations of the sub-probability measure and of the associated SDE equation at the microscale: on the one hand, as a measure-valued solution of a Feynman–Kac-type equation; on the other hand, as the sub-probability associated with an SDE defined up to a survival time with a reaction-dependent rate. These different interpretations give rise to two different microscopic stochastic models and therefore to two different techniques of probabilistic analysis. Finally, by considering the interacting particle systems associated with both models, we discuss how their empirical densities provide two different kernel estimators for the PDE solution.
\end{abstract}

\section{Preliminaries and Motivation}\label{sec:intro}
\emph{Probabilistic interpretations of PDEs} relate solutions of such equations to stochastic processes at a lower scale, given as solution of nonlinear stochastic differential equations (SDEs). Hence, originally deterministic evolution problems can be interpreted through a stochastic perspective.

The relationship between stochastic processes and partial differential equations has played a fundamental role in the development of stochastic analysis \cite{Frejdlin,Friedman,Kac2,Rosenblatt}. In the context of McKean-type advection–diffusion equations, we recall the seminal results due to McKean \cite{McKean1967Propagation}, and, between many others, M\'{e}l\'{e}ard et al. \cite{Meleard_Coppoletta,Meleard} and Sznitman \cite{Snitzman}. Recently this framework has been extended to encompass reaction–advection–diffusion equations \cite{2016_Russo,2019_russo_proceedings,russonumerico,arxivTarquiniUgolini1,2025_morale_tarquini_ugolini}.

To make things precise, we begin by considering the following one-dimensional PDE
\begin{equation}\label{eq:Fokker_Planck}
    \begin{aligned}
        &\partial_t u_t = \frac{\sigma^2}{2}\Delta u_t - \nabla\cdot \big(b(t,x,u_t)u_t\big),\\
        &u_{t=0} = u_0,
    \end{aligned}
\end{equation}
with $(t,x)\in[0,T]\times\mathbb{R}$, where $u_0$ is a given probability density, and $b$ a sufficiently regular function depending on $u$ itself.  Equation \eqref{eq:Fokker_Planck}  is referred to as \emph{McKean-Vlasov-Fokker-Planck equation}.

We start by introducing the general setting. From the probabilistic point of view, the usual interpretation is the following: when the initial condition $u_0$ is a probability density, equation \eqref{eq:Fokker_Planck} describes the evolution of the time marginals of the law of a stochastic process solution to the McKean-Vlasov SDE
\begin{equation}\label{eq:MKV_SDE_intro}
    \begin{aligned}
        &X_t = X_0 + \int_0^t b(s,X_s,\mu_s)ds + \sigma W_t,\\
        &X_0\sim u_0,\\
        &\mu_t:=\mathcal{L}(X_t),
    \end{aligned}
\end{equation}
with $t\in[0,T]$, where $\mathcal{L}(X_t)$ denotes the  probability law of $X_t$. More precisely, given $X$ solution to the McKean-Vlasov SDE \eqref{eq:MKV_SDE_intro}, $t\mapsto\mu_t:=\mathcal{L}(X_t)$ is a measure-valued solution to \eqref{eq:Fokker_Planck}.  Conversely, if $u_t$ is a weak solution of \eqref{eq:Fokker_Planck}, at each time $t\in[0,T]$, the measure $u_t(\cdot)dx$ is the law of $X_t$, solution of \eqref{eq:MKV_SDE_intro}.

For this reason, these time-evolution PDEs are commonly described as \emph{conservative}, meaning that their solutions $u_t $ satisfy the property that
\begin{equation*}
    \int_{\mathbb{R}} u_t(x)dx = 1,\qquad \forall\,t\in[0,T].
\end{equation*}
Within the probabilistic interpretation of McKean-type PDEs, it is common to introduce the interacting particle system associated with the McKean-Vlasov SDE  \eqref{eq:MKV_SDE_intro}. From both a statistical and a numerical viewpoint, the empirical distribution of such a many-particles system is expected to be a asymptotically consistent statistical estimator of the law of the limiting process and therefore provides an approximation of the solution of the corresponding PDE.

Specifically, on the product probability space $\left(\Omega^N,\mathcal{F}^{\otimes N},\mathbb{P}^{\otimes N}\right)$, one can consider the \emph{system of $N\in\mathbb{N}$ mean-field interacting particles} associated with \eqref{eq:MKV_SDE_intro}, defined by
\begin{equation}\label{eq:particle_system_intro}
    X_t^{i,N} = X_0^{i,N} + \int_0^t b(s,X_s^{i,N},\mu_s^N)ds + \sigma W^{i,N}_t,\qquad i=1,\dots,N,
\end{equation}
where  $(W^{1,N},\dots,W^{N,N})$ are independent Brownian motions, and $(X_0^{1,N},\dots,X_0^{N,N})$ are independent random variables which are identically distributed according to the initial law $u_0(\cdot)dx$ and independent of the Brownian motions.

For each $t\in[0,T]$, the interaction between particles is encoded through the \emph{empirical measure}
\begin{equation}\label{eq:intro_emp_measure}
    \mu_t^N := \frac{1}{N} \sum_{i=0}^N \varepsilon_{X_t^{i,N}},
\end{equation}
where $\varepsilon_x$ denotes the Dirac measure at $x\in \mathbb R$.

By construction, the particle system is exchangeable, reflecting the mean-field nature of the interaction. The empirical measure $\mu^N$ plays a central role, as it provides the particle-based approximation of the law $\mu=\mathcal{L}(X)$ of the limiting McKean-Vlasov process. A fundamental question is therefore whether $\mu^N$ converges to $\mu$ as the number of particles $N$ grows to infinity. The classical convergence property considered in this context is the so-called \emph{propagation of chaos} \cite{Snitzman}. More precisely, propagation of chaos holds if
\begin{equation*}
    \mu^N\xrightarrow{L^2} \mu
\end{equation*}
which in turn implies that, for any fixed $k\geq 1$,
\begin{equation*}
    (X^{1,N},\dots X^{k,N})\xrightarrow[N\rightarrow\infty]{d} (X^1,\dots X^k),
\end{equation*}
where $(X^1,\dots,X^k)$ are $k$ independent copies of the limiting process $X$. 
This result formalises the idea that, although particles interact through the empirical measure, they become asymptotically independent as $N\to\infty$, each following the law prescribed by the limit McKean-Vlasov equation \eqref{eq:MKV_SDE_intro}. Hence, a law of large number apply. Moreover, the particle system \eqref{eq:particle_system_intro} can be fully described by such an equation for a representative particle, as all particles behave identically and each ``feels'' the interaction with the infinitely many others as a self-generated field. Finally, from a physical perspective, such particles systems and the related convergence results provide a microscopic interpretation of physical phenomena that are otherwise modeled at a macroscopic level by a nonlinear PDE, such as Equation \eqref{eq:Fokker_Planck}.

Whenever \eqref{eq:Fokker_Planck} includes a reaction term, one obtains a more general class of PDEs than the classical McKean–Fokker–Planck equation previously considered, given by
\begin{equation}\label{eq:Fokker_Planck_reaction}
    \begin{aligned}
        &\partial_t u_t = \frac{\sigma^2}{2}\Delta u_t - \nabla\cdot \big(b(t,x,u_t)u_t\big) + c(t,x,u_t)u_t;\\
        &u_{t=0} = u_0,
    \end{aligned}
\end{equation}
with $(t,x)\in[0,T]\times\mathbb{R}$, where $u_0$ is a given probability measure. 
Such time-evolution PDEs are referred to as \emph{non conservative} in the sense that their solutions $u_t$ no longer verify the property that the probability mass is equal to 1 for each time $t\in[0,T]$.

In the simpler case $b=b(t,x)$ and $c=c(t,x)$, there is the classical result (see, e.g., \cite[Chapter 5, Theorem 7.6]{Karatzas}) stating that the strong solution $v(t,x)$ of the Kolmogorov backward equation
\begin{equation*}
    \begin{aligned}
        &-\partial_t v_t = \frac{\sigma^2}{2}\Delta v_t + \nabla\cdot \big(b(t,x)v_t\big) + c(t,x)v_t,\qquad &(t,x)\in[0,T)\times\mathbb{R};\\
        &v(T,x)=f(x),\qquad &x\in\mathbb{R}
    \end{aligned}
\end{equation*}
admits the probabilistic representation
\begin{equation*}
    v(t,x)=\mathbb{E}\left[f(X_T)\exp\left(\int_t^Tc(s,X_s)ds\right)\Big|X_t = x\right],
\end{equation*}
where $X$ is the solution to the SDE, for any $s\in(t,T]$
\begin{equation*}
    \begin{aligned}
        &dX_s = b(s,X_s)ds + \sqrt{2}dW_s;
        &X_t = x.
    \end{aligned}
\end{equation*}
The extension of the probabilistic representation of McKean-type equations to a large class of non-conservative PDEs has been done only recently in \cite{2017_russo_barbu,2016_Russo,2019_russo_proceedings}.  The authors introduce a representation based on a generalised nonlinear SDE called McKean-Feynman-Kac equation (MKFK-SDE), whose coefficients at time $t\in[0,T]$ do not only depend on the process but also on the solution to a Feynman-Kac-type equation in which a discount rate dependent on the reaction coefficient of the PDE model appears. The reason why such an equation must be solved, rather than having an explicit representation as in the standard Feynman–Kac theory, is that the reaction coefficient depends on the solution itself. On the other hand, the dependence of the drift coefficient on the solution to such an equation reflects the fact that the non-conservative nature of PDE model affects the dynamics of the stochastic equation. Moreover, it is the solution to the above-mentioned Feynman-Kac equation which solves in the distributional sense the PDE, rather than the law of the solution process.

The aim of the present paper is to propose a conceptual viewpoint on the probabilistic interpretation of the PDE as a macroscale deterministic dynamics emerging from stochastic microscale behaviour. Within this perspective, we present and compare two alternative probabilistic approaches in the presence of non-conservative and path-dependent coefficients, as  in \cite{arxivTarquiniUgolini1,2025_morale_tarquini_ugolini}.
We emphasise both their common features and their fundamental differences, with the goal of clarifying the distinct probabilistic interpretations they induce for the same underlying PDE. 
Indeed, the discount factor in the Feynman–Kac interpretation coincides with the survival probability up to time $t$, with a rate determined by the reaction term. From this perspective, a more natural microscopic description is given by a particle that evolves until a random reaction time, modeled as a stopping time $\tau$, at which point it is removed from the system. This consists in introducing a McKean–Vlasov–type SDE, whose solution we denote by $Y$, coupled with a killing time $\tau$. Since the coefficients of the PDE depend on the solution itself, as in the previous setting, the associated SDE and killing time do not involve the law of the state process alone, but rather the joint law of the state and its survival, that is, the distribution of $Y_t$ intersected with the event that the killing time has not yet occurred. In other words, it is the  the sub-probability measure $\nu_t:=\mathbb{P}\left(Y_t\in\cdot,\tau>t\right)$ which solves in the distributional sense the PDE model.

Summarizing, in order to provide a probabilistic interpretation of a McKean-Fokker-Planck equation with reaction, one can rely on two distinct approaches which are based on two different methodologies. The first one could be defined as more analytical, while the second one as more probabilistic. Also, we want to stress here that, when considering the interacting particle systems associated with these two stochastic models, the corresponding empirical densities can be interpreted as two different estimators of the solution to the same PDE.

The target macroscale dynamics that we consider as a case study is given by
\begin{equation*} 
	\begin{split}
		\partial_t  \rho   &=	\Delta\rho - \nabla  \cdot(f(c)  \rho )  -\lambda \rho c,   \\
		\partial_t c & = -\lambda c\, K \ast \rho ,  
	\end{split}
\end{equation*}
subject to the initial conditions $\rho(0,x)= \rho_0(x);\, c(0,x) = c_0$, with $(t,x)\in [0,T]\times\mathbb{R}$, where $T>0$. $K$ is a regularising kernel that here is taken as a Gaussian kernel, even tough it can be generalised to any smooth kernel $K$ satisfying some regularity conditions (see \cite{arxivTarquiniUgolini1,2025_morale_tarquini_ugolini}). The function $f(c)$, explicitly given by \eqref{eq:f(c)}, denotes a velocity field depending on the underlying concentration $c$. By solving the second equation, the system reduces to a non-conservative, path-dependent PDE, with the drift depending upon the measure $\rho(t,\cdot)dx$ for any $t\in [0,T]$.

The two different approaches for a probabilistic interpretation of such system are studied into details in \cite{arxivTarquiniUgolini1,2025_morale_tarquini_ugolini}. In \cite{arxivTarquiniUgolini1}, we extend \cite{2016_Russo} to the path-dependent case, while in \cite{2025_morale_tarquini_ugolini} we detail the analysis for the processes stopped by the reaction. 

The present study is part of a broader applied project aimed at introducing randomness at the into the mathematical modelling framework of marble degradation in the context of cultural heritage. In particular, the focus is on modelling the reaction of calcium carbonate    under the interaction with diffusing sulphur dioxide    \cite{2007_AFNT2_TPM, 2005_GN_NLA}. In this process, sulphur dioxide  diffuses through the stone pores and interact with the surface of such pores. This causes a chemical reaction that produces gypsum, a much more porous stone than marble. Consequently, the material degrades. 

The first deterministic mathematical model is given by the following PDE-ODE system  is given in \cite{2007_AFNT2_TPM, 2004_ADN, 2005_GN_NLA};  it reads, for $(t,x)\in [0,T]\times\mathbb{R}$,  $T>0$ as 
\begin{equation} \label{eq:modello_Natalini} 
	\begin{cases}
		\partial_t  \rho   =  \Delta\rho - \nabla \cdot \left( \frac{\nabla\varphi(c)}{\varphi(c)}\rho \right) - \lambda c \rho.,   \\
		\partial_t c = -\lambda \rho c,  
	\end{cases}
\end{equation}
with  initial conditions
\begin{equation} \label{eq:modello_Natalini_condition1} 
	\rho(0,x)= \rho_0(x); \quad c(0,x) = c_0.
\end{equation}
 In system \eqref{eq:modello_Natalini}, $\rho$ stands for the sulphur dioxide concentration, $c$ for the calcite density, and $\lambda\in \mathbb{R}_{+}$ is the rate of the reaction; the function $\varphi=\varphi(c)$ is the porosity of the material, which is usually considered a linear function of the calcite density $c$ as
\begin{equation}\label{eq:def_porosity}
    \varphi=\varphi(c)=\varphi_0+\varphi_1 c,
\end{equation}
where $\varphi_0\in \mathbb R_+$, $\varphi_1 \in \mathbb R$ are such that $\varphi_0,\varphi_0+\varphi_1 c_0 \in (0,\bar{\varphi})$ with $\bar{\varphi}\in\mathbb{R}_+$.  
 
For the sake of simplicity, in the present work   we consider a constant initial condition for $c$ in \eqref{eq:modello_Natalini_condition1}; even though, the system is well-posed also in the case $c_0(x)\in (0,C_0)$, such that $C_0-c_0 \in H^1(\mathbb{R})$. As a consequence, for any $(t,x)\in[0,T]\times \mathbb R$, one has that $c(t,x)\in (0,C_0)$ \cite{2005_GN_NLA}. Furthermore, for the well-posedness of the system, it is also required that the initial condition $\rho_0(x) \in (0,S_0)$ and that $\rho_0\in L^2(\mathbb R)$ \cite{2005_GN_NLA,2007_GN_CPDE}. Since it is shown that $\rho(t,x) \in (0,S_0)$ for any $(t,x)\in [0,T]\times\mathbb R$, whenever $S_0\le 1$, the concentration  $\rho$  is a map from $[0,T]\times \mathbb{R}$ to $[0,1]$. This is the case we consider here.

Given the explicit solution of the equation for $c$ as a function of $\rho$:
\begin{equation}\label{eq:c_explicit}
    c(t,x) = c_0\exp\left(-\lambda\int_0^t\rho(u,x)du \right)= c_0\exp\big( -\lambda  \rho(\cdot,x)(t) \big),
\end{equation}
the velocity field 
\begin{equation}\label{eq:f(c)}
f(c)= \frac{\nabla\varphi(c)}{\varphi(c)} 
\end{equation}
exerted on the concentration $\rho$   becomes 
\begin{equation}\label{eq:def_b_PDE}
    b\big(\rho(\cdot,x)(t), \nabla \rho(\cdot,x)(t) \big) := -\varphi_1\lambda c_0\frac{\exp\big( -\lambda \rho(\cdot,x)(t)\big)\nabla \rho(\cdot,x)(t)}{\varphi_0+\varphi_1c_0 \exp\big( -\lambda \rho(\cdot,x)(t)\big)},
\end{equation}
where  $\rho(\cdot,x)(t)$ denotes the integral function of the concentration $\rho$, i.e. for any $(t,x)\in [0,T]\times\mathbb{R}$, 
\begin{equation*}
    \rho(\cdot,x)(t):=\int_0^t \rho(s,x)ds.
\end{equation*}

Therefore, system \eqref{eq:modello_Natalini} is equivalent to the following non conservative McKean-type evolution equation for $\rho$ 
\begin{equation}\label{eq:PDE_rho_only}
        \partial_t\rho (t,x)  = \Delta\rho(t,x)   - \nabla \cdot \left( b\big(\rho(\cdot,x)(t), \nabla \rho(\cdot,x)(t) \big)\rho(t,x) \right) - \lambda c_0\exp\big( -\lambda  \rho(\cdot,x)(t) \big) \rho(t,x).
\end{equation}
  
To the best of our knowledge, degradation phenomena have been modeled either through statistical approaches \cite{2018_saba} or deterministic PDE frameworks (\cite{2007_AFNT_TPM,2007_AFNT2_TPM,2023_Bonetti_natalini_NLA,2019_BCFGN_CPAA,2005_GN_NLA,2007_GN_CPDE}). Recent works introduced stochastic formulations \cite{MACH2021_AGMMU,2025_MauMorUgo,2025_Mach2023_MRU}, including a random PDE version of \eqref{eq:modello_Natalini} with stochastic boundary conditions \cite{2025_MACH2023_AMMU,2025_MauMorUgo}, whose behaviour and discretizations were analysed in \cite{2024_ArceciMoraleUgolini_arxiv,2025_ADMU_arxiv}. At the nanoscale, a first stochastic interacting particle model has been proposed in \cite{2025_JMMRU_arxiv,2025_Mach2023_MRU}.

As explained above, the idea carried out in \cite{arxivTarquiniUgolini1,2025_morale_tarquini_ugolini} is to provide a deeper understanding of the PDE system \eqref{eq:modello_Natalini} as a mean field approximation of a system of Brownian particles subject to a non-conservative evolution law.

Actually, since we carry out the analysis in the case of a McKean-Vlasov-type drift, it is necessary to regularise the drift of the associated SDE by introducing a non-local dependence upon $\rho$. By regularising at the level of the ODE for $c$,
\begin{equation*}
    \partial_t c = -\lambda   (K*\rho) c,
\end{equation*}
where $*$ stands for the convolution with a smooth kernel  $K: \mathbb{R} \rightarrow \mathbb{R}$, Equation \eqref{eq:PDE_rho_only} becomes
\begin{equation}\label{eq:PDE_rho_regularised}
    \begin{aligned}
        \partial_t\rho(t,x) \,=\, &\Delta\rho(t,x) - \nabla \cdot \left[ b\big(K* \rho(\cdot,x)(t), \nabla K* \rho(\cdot,x)(t)\big) \rho\, \right]\\
        &-\lambda c_0\exp\big( -\lambda  K*\rho(\cdot,x)(t) \big) \rho(t,x).
    \end{aligned}
\end{equation}
\begin{remark}
    At least formally,  in the case in which the kernel $K$ is   a Dirac delta function, the non-local equation \eqref{eq:PDE_rho_regularised} becomes the local one \eqref{eq:PDE_rho_only}.
\end{remark}

The paper is organised as follows. In Section \ref{sec:SDE_PDE_link}, we introduce the stochastic equations associated to the PDE \eqref{eq:PDE_rho_regularised},\eqref{eq:def_b_PDE}.  The particle systems related to such SDEs are introduced in Section \ref{sec:associted_PS}, and the peculiar fact that they essentially represent to different estimators of the solution to the PDE model is pointed out. Finally, in Section \ref{sec:proofs_well_posed}, we briefly outline the main ideas of the proofs of existence and pathwise uniqueness of the strong solution of the two stochastic models introduced in Section \ref{sec:SDE_PDE_link}.

\paragraph{Notation.} Let $(\Omega, \mathcal{F},\mathbb{F}=\{\mathcal{F}_t\}_{t\in [0,T]},\mathbb{P})$ be a filtrated probability space and choose $T \in \mathbb R_+$ to be a sufficiently large time horizon. All stochastic processes throughout the paper are assumed to be defined and adapted to $\mathbb{F}$. Moreover, for any stochastic process $X$, $\mathcal{L}(X)$ denotes its law, and $\mathcal{L}(X_t)$ denotes the time-marginal law of $X$ at time $t\in[0,T]$. 

For any Banach space $V$, we denote with $\mathcal{B}V$ the Borel $\sigma$-algebra on V.

We indicate with $\{\Delta\}$ a \emph{cemetery state} isolated from $\mathbb{R}$ and we consider the convention that any function $f$ defined on $\mathbb{R}$ can be extended to $\mathbb{R}\cup\{\Delta\}$ by setting $f(\Delta)=0$.

\section{ Stochastic equations at the microscale}\label{sec:SDE_PDE_link}
Here we present and discuss two stochastic models associated with the PDE \eqref{eq:PDE_rho_regularised},  with drift \eqref{eq:def_b_PDE}, emphasizing the alternative probabilistic treatments of the reaction term.

In \cite{arxivTarquiniUgolini1,2025_morale_tarquini_ugolini}, we consider the following SDE
\begin{equation}\label{eq:SDE_gamma2}
        Y_t = Y_0 + \int_0^t b\left(\int_0^s(K*\nu_r)(Y_s)dr, \int_0^s(\nabla K*\nu_r)(Y_s)dr\right)ds + \sqrt{2}W_t,
\end{equation}
with $t\in[0,T]$, subject to the initial condition $Y_0\sim\rho_0(x)dx$, where $\nu_t$ is a sub-probability measure. The reason why a sub-probability measure appears in the dynamics of the associated SDE model is the fact that the PDE is of non-conservative McKean type: the coefficients, including the reaction one, depend on the solution of the equation itself. In other words, the non-conservative nature of the PDE carries over to the associated SDE dynamics, and the measure $\nu_t$ represents the law of $Y_t$ up survival at time $t$.

In the following, we discuss how to characterize the measure $\nu_t$. Following a Feynman-Kac approach, it is possible to regard $\nu$ as the solution to a measure-valued Feynman-Kac-type equation. This is done in \cite{arxivTarquiniUgolini1} and explained in the next subsection. On the other hand,  $\nu$ can be modeled as the joint law of the state of the process $Y$ \eqref{eq:SDE_gamma2} and its survival with respect to a killing mechanism. This is done in \cite{2025_morale_tarquini_ugolini} and the main ideas are reported in Subsection \ref{subsec:second_approach}.

\subsection{A Feynman-Kac approach}
The sub-probability measure $\nu$ is considered here as the measure-valued solution to the following Feynman-Kac-type equation
\begin{equation}\label{eq:def_gamma_m2}
        \nu_t(f)  = \mathbb E\left[f\left( Y_t\right)\exp\left(-\lambda c_0 \int_0^t \exp\left( -\lambda\int_0^s(K*\nu_r)(Y_s)dr\right)ds\right)\right],
\end{equation}
for every $t\in[0,T]$ and for every $f\in C_b(\mathbb{R})$.

The fully stochastic system \eqref{eq:SDE_gamma2}-\eqref{eq:def_gamma_m2} admits a pathwise unique strong solution \cite{arxivTarquiniUgolini1}. The key arguments of the proof may be found in Section \ref{sec:proofs_well_posed}. 
 
 A reduction in the complexity of the system \eqref{eq:SDE_gamma2}–\eqref{eq:def_gamma_m2} is achieved by observing that the measure  $\nu_t$
 coincides with the solution to the following measure-valued PDE
\begin{equation}\label{eq:PDE_gamma^m}
 \partial_t\nu_t  = \Delta\nu_t -\nabla\cdot\left[ b\left(\int_0^t u_K(s,\cdot) ds,\int_0^t \nabla u_K(s,\cdot) ds\right)\nu_t\right]-\lambda c_0\exp\left( -\lambda \int_0^t  u_K(s,\cdot) ds \right)\nu_t,
\end{equation}
with $\nu_0(dx)=\rho_0(x)dx$, for $(t,x)\in[0,T]\times\mathbb{R}$, where
\begin{equation}\label{eq:def_mollified_density}
    u_K(t,x)=K*\nu_t(x),\qquad \forall\,(t,x)\in [0,T]\times \mathbb R
\end{equation}
is called the \emph{mollified density} associated to $\nu$.

Equation \eqref{eq:PDE_gamma^m} can be read as a measure-valued version of the deterministic model \eqref{eq:PDE_rho_regularised}. Therefore, it is no longer the law of the process $Y$ but rather the solution of equation \eqref{eq:def_gamma_m2}, in which the law of $Y$ is discounted at a rate given by the reaction coefficient, which is associated to the PDE model. Moreover, the solution to the PDE admits a probabilistic characterisation in terms of a Feynman–Kac–type equation, rather than an explicit representation as in the classical Feynman–Kac formula, since the reaction rate depends on the solution itself. Similarly, the dependence of the advection term in the PDE on the solution implies that the stochastic dynamics, through the drift in \eqref{eq:SDE_gamma2}, are influenced by the loss of mass induced by the reaction term.

The above probabilistic interpretation of the PDE model \eqref{eq:PDE_rho_regularised} suggests a heuristic physical interpretation by means of the equations \eqref{eq:SDE_gamma2}-\eqref{eq:def_gamma_m2}, illustrating how such an SDE model and the deterministic equation \eqref{eq:PDE_rho_regularised} describe the sulphation reaction at different scales.

From a physical perspective, for any $t\in [0,T]$, the measure  $\nu_t$ is a macroscopic observable representing the total amount of sulphur dioxide at time $t$, while $u_K(t,\cdot)$ describes a non-local average of this quantity obtained through convolution with $K$, thereby introducing further non-local effects. In particular, $u_K$  depends on the distribution of the measure over the entire domain, rather than solely on the local behaviour around any point $x\in \mathbb R$.

At the microscale, the stochastic process $Y$ whose dynamics can be rewritten for any $t\in [0,T]$ as 
\begin{equation}\label{eq:SDE_UK_FK} 
    Y_t = Y_0 + \int_0^t b\left(u_K(\cdot,Y_s)(s),\nabla  u_K(\cdot,Y_s)(s)\right)ds + \sqrt{2}W_t
\end{equation}
evolves under the influence of the non-local, discounted, and regularized field $u_K$, given by \eqref{eq:def_mollified_density}, acting at the macroscale. Note that from \eqref{eq:def_gamma_m2}  with $f(\cdot)=K(y-\cdot)$ and \eqref{eq:def_mollified_density}, the mollified density may be written as
\begin{equation}\label{eq:mollified_density_FK}
u_K(t,y)  = \mathbb E\left[K(y-Y_t)\exp\left(-\lambda c_0 \int_0^t \exp\left( -\lambda\int_0^su_K(r,Y_s)dr\right)ds\right)\right],
\end{equation}
Hence, the SDE \eqref{eq:SDE_gamma2} (or, equivalently, \eqref{eq:SDE_UK_FK}) may be coupled at any time $t\in[0,T]$ with the deterministic solution of the PDE \eqref{eq:PDE_gamma^m}, reducing the
complexity of the full system. Alternatively, it could be coupled with the solution of the \emph{measure-valued McKean-Feynman--Kac equation} \eqref{eq:def_gamma_m2} (or its mollified version \eqref{eq:mollified_density_FK}), which represents the underlying field averaged at the macroscopic scale given the process $Y$ at the microscopic level.

The process $Y$ solution to \eqref{eq:SDE_gamma2} can be regarded as the dynamics of a typical particle at the microscale.  The fully stochastic process is the so called  \emph{  McKean-Feynman-Kac SDE} (MKFK-SDE), given by \eqref{eq:SDE_gamma2} coupled with  \eqref{eq:def_gamma_m2}; while the hybrid deterministic-stochastic model is given by \eqref{eq:SDE_gamma2} coupled with the deterministic equation \eqref{eq:PDE_gamma^m}.
  
As it is well known, the discount rate appearing in \eqref{eq:def_gamma_m2} corresponds to the survival probability of a killing mechanism. In fact, let us consider the following stopping time $\tau$  
\begin{equation}\label{eq:stopping_time_kill}
    \tau \,:=\, \inf_{t\geq 0}\Bigg\{ \int_0^t \lambda c_0\exp\left(-\lambda\int_0^s K*\nu_r(Y_s)dr \right)ds \,\geq\, Z \Bigg\},     
\end{equation}
with $Z\sim Exp(1)$ independent of $Y_0$ and $W$. Then, the measure $\nu_t$ satisfies the following condition
\begin{equation}\label{eq:def_subprobability_nu_t}
   \nu_t = \mathbb{P}\left( Y_t \in \cdot,\ t < \tau \right).
\end{equation}
Indeed, for any test function $f\in C_b^{\infty}(\mathbb{R})$,
\begin{equation*}
    \begin{aligned}
        \int_\mathbb{R} f(x)\mathbb{P}\left( Y_t \in dx,\ t < \tau \right) &=\mathbb{E}\left[f(Y_t) \exp\left(-\int_0^t \lambda c_0\exp\left( -\lambda \int_0^s K*\nu_r(Y_s)dr \right)ds \right)\right]\\
        &= \int_\mathbb{R} f(x)\nu(dx),
    \end{aligned}
\end{equation*}
where the second equality is due to \eqref{eq:def_gamma_m2}.

In view of the interpretation of the discount rate in the Feynman–Kac–type equation as a killing rate, it is natural to seek a probabilistic interpretation of the deterministic model via the McKean–Vlasov SDE   \eqref{eq:SDE_gamma2} with killing \eqref{eq:stopping_time_kill}.  

\subsection{A killed McKean-Vlasov SDE}\label{subsec:second_approach}
Starting from \eqref{eq:def_subprobability_nu_t}, a more natural microscopic description consists of a particle evolving according to the dynamics \eqref{eq:SDE_gamma2} up to the random stopping time $\tau$ \eqref{eq:stopping_time_kill}, after which it is removed from the system and sent to a cemetery state $\{\Delta\}$. 

The microscopic dynamics $Y$ is then described by the following system
\begin{equation}\label{eq:model_killing_measure}
    \begin{aligned}
        &Y_t = Y_0 + \int_0^t b\left(\int_0^s K*\nu_r(Y_s)\,dr,\ \int_0^s \nabla K*\nu_r(Y_s)\,dr\right)\,ds + \sqrt{2} W_t, \quad &t < \tau;\\
        &Y_t\in\{\Delta\};\quad &t\geq \tau,
    \end{aligned}
\end{equation}
 where $\nu_t$  is the law \eqref{eq:def_subprobability_nu_t} of $Y$ up to survival  and $\tau$ is given by \eqref{eq:stopping_time_kill}.
 Also, it is assumed the convention $Y_{\tau}:=\lim_{t\rightarrow\tau-}Y_t$.
The stochastic model \eqref{eq:model_killing_measure} admits a pathwise unique strong solution \cite{2025_morale_tarquini_ugolini}.  See Section \ref{sec:proofs_well_posed} for the main ideas.

 System \eqref{eq:model_killing_measure} is fully stochastic and it gives a description at the microscale of the dynamics. Again one may show that, for any $t\in [0,T]$, the measure $\nu_t$ ( and not the law of $Y_t$)  is a solution to the measure-valued version \eqref{eq:PDE_gamma^m} of the deterministic model \eqref{eq:PDE_rho_regularised}.  By standard arguments, as the Hahn-Banach theorem and Riesz's representations theorems \cite{Meleard_Coppoletta,2025_morale_tarquini_ugolini}, one may prove that the measure $\nu_t$ is regular enough to admits a density  $v(t,\cdot)$. Furthermore, for any  $f\in C_b^2\left(\mathbb{R}\right)$ and \(t\in[0, T]\) we obtain
\begin{align*}
    \int_{\mathbb{R}}f(x)v(t,x)dx = &\int_{\mathbb{R}}f(x)\rho_0(x)dx \,+\, \int_0^t\int_{\mathbb{R}}\Delta f(x)v(s,x)dxds&&\\
    &+\int_0^t\int_{\mathbb{R}}\nabla f(x)b\left(K*v(\cdot,x)(s),\nabla K*v(\cdot,x)(s)\right)v(s,x)dxds&&\\
    &-\int_0^t \int_{\mathbb{R}}  f(x)\lambda c_0\exp\big( -\lambda K*v(\cdot,x)(s)\big)v(s,x)dxds.
\end{align*}
Therefore, $\{v(t,\cdot)\}_{t\in[0,T]}$ is a weak solution to the equation \eqref{eq:PDE_rho_regularised}, \cite{2025_morale_tarquini_ugolini}.

Again, the dynamics of the system under consideration may be seen from two complementary viewpoints. On the one hand, a fully microscopic description is provided by a purely stochastic dynamics, governed by the stopped SDE
\begin{align}
    Y_t &= Y_0 + \int_0^t b\left(\int_0^s K*v(r,Y_s)\,dr,\int_0^s \nabla K*v(r,Y_s)\,dr\right)\,ds
    + \sqrt{2}\,W_t,\qquad t<\tau,\label{eq:SDE_kill_density}\\
    v(t,\cdot)\,dx &= \mathbb{P}(Y_t\in dx,\ t<\tau)=:\nu_t(dx),\label{eq:SDE_kill_density_law}
\end{align}
where the killing time $\tau$ is defined in \eqref{eq:model_killing_measure}.

On the other hand, one may consider a hybrid description in which the stochastic dynamics of the particle is coupled with a deterministic evolution equation for the density $v$, namely
\begin{align}
   Y_t &= Y_0 + \int_0^t b\left(\int_0^s K*v(r,Y_s)\,dr,\ \int_0^s \nabla K*v(r,Y_s)\,dr\right)\,ds
   + \sqrt{2}\,W_t,\qquad t<\tau,\label{eq:SDE_kill_density2}\\
   \partial_t v(t,x)
   &= \Delta v(t,x)
      - \nabla \cdot \Big[ b\big(K* v(\cdot,x)(t), \nabla K* v(\cdot,x)(t)\big)\, v(t,x) \Big] \notag\\
   &\quad - \lambda c_0 \exp\!\big( -\lambda K* v(\cdot,x)(t) \big)\, v(t,x),
   \label{eq:SDE_kill_density2_bis}
\end{align}
for $t\in[0,T]$, where \eqref{eq:SDE_kill_density2_bis} is understood in the distributional sense.
\medskip

We conclude by providing a heuristic physical interpretation of the SDE model \eqref{eq:SDE_kill_density}.
Here, the dependence of the dynamics of the process $Y$ on the joint law of $Y$ itself and its survival accounts for the mean-field interaction of every sulphur dioxide  molecule with the others that have not reacted yet. The stopping time $\tau$ takes care of the interaction of a  molecule with a calcium carbonate one; when such a reaction happens, a gypsum molecule  is formed and the killing term reflects the fact that the sulphur dioxide particle is consumed by the chemical reaction.

\section{Particle approximations of the PDE solution}\label{sec:associted_PS}

We now consider a particle system (PS) associated with the stochastic microscale dynamics and discuss how the two approaches introduced in the previous section provide two possible consistent estimators of the solution to the PDE \eqref{eq:PDE_rho_regularised}, or, equivalently, two distinct estimators of the sub-probability measure  $\nu_t$ in \eqref{eq:def_subprobability_nu_t}. In order to avoid discussing regularity issues of the interaction kernel, we may take $K$ be the Gaussian kernel. See \cite{arxivTarquiniUgolini1,2025_morale_tarquini_ugolini} for a generic smooth kernel $K$ satisfying suitable assumptions.
 
Let $N\in\mathbb{N}$ and let us introduce the product probability space $\left(\Omega^N,\mathcal{F}^{\otimes N},\mathbb{P}^{\otimes N}\right)$, endowed with the filtration obtained as the natural extension of the original filtration to the product space. On this probability space, let $(W^{1,N},\dots,W^{N,N})$ be a family of $N$ independent $\mathbb{F}$-Brownian motions, and let $(Y^{1,N},\dots,Y^{N,N})$ denote an $\mathbb{R}^N$-valued stochastic process representing the positions of $N$ particles. The initial conditions are given by $Y_0^{i,N}=Y_0^i$, where $(Y_0^1,\dots,Y_0^N)$ are independent random variables, identically distributed as $\rho_0(x)\,dx$, which are independent of the Brownian motions $\{W^{j,N}\}_{j=1}^N$.

Each particles undergoes to the following dynamics 
\begin{equation} \label{eq:particle-dynamics}
  Y_t^{i,N}  = Y_0^{i,N} + \int_0^t b\left(u_N\left( \cdot,Y_s^{i,N} \right)(s),\nabla u_N\left( \cdot,Y_s^{i,N} \right)(s)\right)ds + \sqrt{2}\,W_t^i, \qquad  t\in I^i.\end{equation}
The different approaches introduced in the previous section specify the time interval $I^i$ on which the equation is defined for the 
$i-$th particle in the system, as well as the empirical density $u_N$.
\smallskip

\emph{The MKFK-PS.} Let $I^i=[0,T]$ for any particle $i$. A natural estimator of the mollified density  \eqref{eq:mollified_density_FK} is then the following 
\begin{equation}\label{eq:system_particles_interacting_density}
    u_N(t,y) =\frac{1}{N}\sum_{i=1}^N K\left( y - Y_t^{i,N} \right) \exp\left(-\lambda c_0 \int_0^t \exp\left(-\lambda \int_0^s u_N\left( r,Y^{i,N}_s \right) dr\right) ds\right),
\end{equation}
Consistently, the above empirical mollified density is a weighted convolution of the kernel with the empirical particle measure, that is, each term of the discrete convolution is weighted by the discount term.  

The convergence of $u_N$ to $u_K$ (see \cite[Proposition 4.4]{arxivTarquiniUgolini1}), together with the probabilistic interpretation discussed in Section \ref{sec:SDE_PDE_link}, implies that $u_N$ provides a statistical estimator of the solution to $K*\rho$ where $\rho$ is the solution to the PDE model with the regularizing kernel $K$ appearing in the coefficients. We emphasise that such convergence cannot be expected without such a smoothing. Moreover, this convergence result also entails propagation of chaos for the associated interacting particle system \cite[Theorem 4.5]{arxivTarquiniUgolini1}.

Precisely, let us consider $N\in\mathbb{N}$ independent copies of the SDE \eqref{eq:SDE_UK_FK}-\eqref{eq:mollified_density_FK}. That is,
\begin{equation}\label{eq:system_particles_iid}
    \begin{aligned}
        Y_t^i &= Y_0^i + \int_0^t b\big(u_K\left(\cdot,Y_s^i\right)(s),\nabla u_K\left(\cdot,Y_s^i\right)(s)\big)ds + \sqrt{2}W_t^i;\\
        u_K(t,y) &= \mathbb{E}\left[ K\left(y - Y_t^i\right)\exp\left(-\lambda\,c_0 \int_0^t \exp\left( -\lambda u_K(\cdot,Y^i_s)(s)  \right) ds\right) \right],
    \end{aligned}
\end{equation}
for $t\in[0,T]$ and $i=1,\dots,N$.

For any given  $N\in \mathbb N$, let $(Y^1,\ldots,Y^N)$  be the solution of  system \eqref{eq:system_particles_iid} and  $(Y^{1,N},\ldots,Y^{N,N})$ be the one of \eqref{eq:particle-dynamics}-\eqref{eq:system_particles_interacting_density} with $I^i=[0,T], i=1,...,N$,  respectively. It holds that, for any $k\ge 2$,
\begin{equation*}
    \left( Y^{1,N},Y^{2,N},\ldots,Y^{k,N} \right)\xrightarrow[\mathcal{L}]{N\to \infty} \left( Y^1,Y^2,\ldots,Y^k \right).
\end{equation*}
This means that, in the infinite-particle limit, the full interacting particle system converges to independent particle dynamics, allowing the system to be completely described by the associated SDE. From a physical perspective, this justifies the interpretation of the process $Y$ solution to \eqref{eq:SDE_gamma2}  as the dynamics of a typical \ch{SO_2} molecule at the microscale.

Equivalently, it is well-known that the propagation of chaos property can be rewritten as
\begin{equation*}
    \mathcal{L}\left( Y^{1,N},Y^{2,N},\ldots,Y^{k,N} \right)\xrightarrow{N\to \infty}\otimes_k m,
\end{equation*}
where $m:=\mathcal{L}(Y)$ with $Y$ the solution to SDE \eqref{eq:SDE_gamma2}. That is, the joint law of any $k$-th dimensional process solution of \eqref{eq:particle-dynamics}-\eqref{eq:system_particles_interacting_density} converges to the $k$-dimensional product law of $k$ independent copies of the law $m$ of the solution of \eqref{eq:system_particles_iid}.
\medskip
 
\emph{The PS with killing.}  Let $I^i:=[0,\tau^i)\cap [0, T]$, where $\tau^i$ is an $\mathbb{F}$-stopping time, and $i=1,\dots,N$.  Hence each particle follows the dynamics \eqref{eq:particle-dynamics} up to time $\tau^i$.

Let $\Gamma^N_t\subseteq\{1,\dots,N\}$ the subset of indexes corresponding to alive particles  out of $N\in \mathbb N$ at time $t$. A natural estimator of the sub-probability measure \eqref{eq:def_subprobability_nu_t} is given by  the empirical measure  
\begin{eqnarray}\label{eq:empirical_measure_mu}
\nu_t^N (\cdot):=\frac{1}{N}\sum_{i\in\Gamma_t^N} \varepsilon_{Y^{i,N}_t}(\cdot)=\frac{1}{N}\sum_{i=1}^N \varepsilon_{Y^{i,N}_t}(\cdot)\mathbbm{1}_{(t,\infty)}(\tau^i)\end{eqnarray}
 that is the empirical measure for the surviving particles. As a consequence, the empirical density is
\begin{equation}\label{eq:empirical_density_killed_PS}
    u_N\left(t,x\right):= \left(K*\nu_t^N\right)(x) := \int_{\mathbb{R}}K(x-y)\nu_t^N(dy)=\frac{1}{N}\sum_{i=1}^NK\left(x-Y^{i,N}_t\right)\mathbbm{1}_{(t,\infty)}(\tau^i).
\end{equation}
The stopping time $\tau^i$ is the lifetime of the $i$-th particle and it is defined as
\begin{equation}\label{eq:stopping_time_PS}
    \tau^i = \tau^i\left(Y^{i,N}\right) \,:=\, \inf_{t\geq 0}\Bigg\{ \int_0^t \lambda c_0\exp\left(-\lambda \int_0^s u_N(r,Y^{i,N}_s)dr \right)ds \,\geq\, Z_i \Bigg\},
\end{equation}
where $Z_i\sim Exp(1)$ is independent of $Y^{1,N},\dots,Y^{N,N}$, for every $i=1,\dots,N$.

Well-posedness of system \eqref{eq:particle-dynamics},\eqref{eq:system_particles_interacting_density} and  \eqref{eq:particle-dynamics},\eqref{eq:empirical_density_killed_PS}, has been established in \cite{arxivTarquiniUgolini1} and \cite{2025_morale_tarquini_ugolini}, respectively. Heuristically, it is easy to see that if a limit measure $\mu_t(dx)= v(t,\cdot)dx$  as $N$ increasing to infinity exists, then $v(t,\cdot) $ is a weak solution of \eqref{eq:PDE_rho_regularised}. Rigorous derivation is addressed in a future work.

\smallskip
In conclusion, for a given McKean-type PDE with reaction, the microscale dynamics can be interpreted in different ways. By looking at the empirical densities \eqref{eq:system_particles_interacting_density} and \eqref{eq:empirical_density_killed_PS} related to the two particle systems, we observe that they are different in one key way. In the first one, all the particles are alive for every time $t\in[0,T]$ and the reaction in the PDE model is accounted for by introducing a discount factor that depends on the reaction coefficient of the deterministic model itself; in the second one, at every time $t\in[0,T]$, if the killing time $\tau^i$ (whose intensity is the reaction coefficient of the PDE model) relative to the $i$-th particle has stricken, such a particle is removed. Therefore, two different estimators are obtained for the solution to the PDE model \eqref{eq:PDE_rho_regularised} convoluted with a regularising kernel $K$. This statistical theoretical procedure is referred to as in the literature as \emph{kernel estimation} \cite{kerneldensityref}.
\section{A note on the well-posedness of the equations}\label{sec:proofs_well_posed}
For completeness, we conclude by discussing the existence and pathwise uniqueness of strong solutions for the two SDE models introduced in Section \ref{sec:SDE_PDE_link}.

 \emph{The MKFK-SDE.} As the  stochastic system \eqref{eq:SDE_gamma2}-\eqref{eq:def_gamma_m2} regards, we note that the expectation in \eqref{eq:def_gamma_m2} can be rewritten as an integral with respect to the path-space law  $m=\mathcal{L}(Y)$. Therefore, \eqref{eq:def_gamma_m2} can be generalised by considering a generic measure $m\in\mathcal{P}(C)$, with $C:=C([0,T],\mathbb{R})$, which is not necessarily the law of $Y$. Hence, for any $t\in [0,T]$, the generalised \emph{McKean-Feynman-Kac}-equation is
 \begin{align}
    &u_{K,m}(t,y) = \int_{C} K(y - X_t(\omega))\exp\left(-\lambda\,c_0 \int_0^t \exp\left(-\lambda \int_0^s u_{K,m}\left( r,X_s(\omega) \right) dr\right) ds\right)m(d\omega)\label{eq:first_model_SDE_2},
\end{align}
subject to the initial condition $ Y_0\sim \rho_0(x)dx$, where $X$ denotes the canonical process.

In \cite[Subsection 2.1]{arxivTarquiniUgolini1}, it is established that \eqref{eq:first_model_SDE_2} admits a unique solution by a fixed point argument. To make this more precise, let $m\in \mathcal{P}(C)$ and let $X:=X^m$ be the associated canonical process. The equation can be naturally interpreted as a fixed point problem in a suitable Banach space of real-valued continuous processes on $[0,T]$. This functional framework allows us to reformulate the problem in an abstract but convenient way. In fact, equation \eqref{eq:first_model_SDE_2} can be rewritten as
\begin{equation}\label{eq:linking_fixed_point}
    u_{K,m} = \left(    T^m \circ\tau\right)\left(u_{K,m}\right),
\end{equation}
where the operator $T^m$ maps a process to a deterministic function obtained by averaging, with respect to the measure $m$, a suitable functional of the process in which the reaction rate appears, while $\tau$ evaluates such a function along the trajectories of the canonical process. In this way, the problem of well-posedness of the Feynman-Kac-type equation \eqref{eq:first_model_SDE_2} is reduced to showing that the composition $\tau\circ T^m$ is a contraction. Once this property is established, the existence and uniqueness of a fixed point follow directly from Banach’s fixed point theorem, yielding the existence of a solution to \eqref{eq:first_model_SDE_2} \cite[Proposition 2.7]{arxivTarquiniUgolini1}. Uniqueness can be proved by standard arguments, along the same lines as in \cite[Theorem 3.1]{2016_Russo}.

 Given the well-posedness of \eqref{eq:first_model_SDE_2}, it is proved that the following SDE
\begin{align}
    &Y_t = Y_0 + \int_0^t b\big(u_{K,m}(\cdot,Y_s)(s),\nabla u_{K,m}(\cdot,Y_s)(s)\big)ds + \sqrt{2}W_t; \quad t\in[0,T]\label{eq:first_model_SDE_1}
\end{align}
admits a pathwise unique strong solution \cite[Proposition 3.3]{arxivTarquiniUgolini1}.
 
Thanks to the boundedness and Lipschitz continuity of the solution to the McKean-Feynman-Kac equation \eqref{eq:linking_fixed_point} and the explicit form of the drift coefficient $b$ \eqref{eq:def_b_PDE}, this can be achieved with the usual arguments to established well-posedness of McKean-Vlasov-type SDEs. Specifically, let $\Theta$ be the map that sends square integrable probability measure over the space of continuous paths $\mathcal{P}^2\left(C \right)$ into the law of the process solution to the SDE \eqref{eq:first_model_SDE_1}-\eqref{eq:first_model_SDE_2}. Fixed a measure $m\in \mathcal{P}^2\left(C \right)$, by iteratively applying the map $\Theta$  to the measure $m$, we define a sequence of measures $  \{m_k\}_{k\in \mathbb{N}}$ in $\mathcal{P}^2\left(C \right)$  such that
\begin{equation}\label{eq:def_m_k_ex_uni_SDE}
    m_k\,:=\,  \Theta^k(m) = \Theta\left(\Theta^{k-1}(m)  \right).
\end{equation} 
    
It is proved that there exists a probability measure $\bar{\mu}\in\mathcal{P}^2\left(C \right)$ such that $m_k$ weakly converges to  $\bar{\mu}$ as $k$ increases to $\infty$, and from \eqref{eq:def_m_k_ex_uni_SDE} we get $\bar{\mu}=\Theta(\bar{\mu})=\mathcal{L}(Y^{\bar{\mu}})$; hence,   the existence of a weak solution to MKFK-SDE \eqref{eq:first_model_SDE_1} coupled with \eqref{eq:first_model_SDE_2} is achieved, and then since it results pathwise unique,   the solution exists in the strong sense (see, e.g.,\cite[pg. 308]{Karatzas}).
\medskip

\emph{The killed SDE.} The technique for the well-posedness of the equation  \eqref{eq:model_killing_measure} with killing time \eqref{eq:stopping_time_kill} is different from the regular case. The main issue being that in \eqref{eq:model_killing_measure} it appears a sub-probability measure instead of a probability law.

The trick we use is to reformulate \eqref{eq:model_killing_measure} by lifting it to an \(\mathbb{R}^2\)-valued process $(\widetilde{X},\widetilde{\Lambda})$, where the first component is the state and the second is the intensity of the killing rate. Specifically, for every $t\in[0,T]$,
\begin{equation}\label{eq:SDE_con_killing_solo_dentro}
    \begin{aligned}
        &\widetilde{X}_t = X_0 + \int_0^t b\left(\int_0^s K*\nu_r(\widetilde{X}_s)dr, \int_0^s\nabla K*\nu_r(\widetilde{X}_s)dr\right)ds + \sqrt{2} W_t;\\
        &\widetilde{\Lambda}_t = \int_0^t \lambda c_0\exp\left( -\lambda\int_0^s K*\nu_r(\widetilde{X}_s)dr\right) ds,
    \end{aligned}
\end{equation}
and $\nu_t=\mathbb P\left( \widetilde{X}_t\in\cdot, \widetilde{\Lambda}_t<Z\right)$.

Then, \eqref{eq:SDE_con_killing_solo_dentro} is rewritten using a functional mapping probability measures on \(\mathbb{R}^2\) to sub-probability measures on \(\mathbb{R}\) (see, e.g., \cite[Proposition 3.12]{2025_hambly}). Therefore, we have obtained a McKean-Vlasov-type SDE in which a true probability measure appears. Thus, the usual arguments can be used to prove strong existence of a pathwise unique solution \cite[Theorem 2.8]{2025_morale_tarquini_ugolini}.
 
Notice that in \eqref{eq:SDE_con_killing_solo_dentro} the dynamics is defined for all times in $[0,T]$. The strong existence and pathwise uniqueness of a solution to the system \eqref{eq:stopping_time_kill} with killing time \eqref{eq:model_killing_measure} is an immediate consequence.

\section*{Acknowledgments}
 D.M. and S.U. are members of GNAMPA (Gruppo Nazionale per l’Analisi Matematica, la Probabilità e le loro Applicazioni) of the Italian Istituto Nazionale di Alta Matematica (INdAM).

\printbibliography

\end{document}